\documentclass[11pt]{amsart}

\usepackage{amssymb}
\usepackage{forloop}
\usepackage{pst-all}
\usepackage{auto-pst-pdf}
\usepackage{calc}

\evensidemargin\oddsidemargin

\begin{document}
\baselineskip=18pt
\setcounter{page}{1}
    
\newtheorem{Corr}{Corollary}
\newtheorem{Prop}{Proposition}
\newtheorem{Theo}{Theorem}
\newtheorem{Lemm}{Lemma\!\!}
\newtheorem{Rq}{Remark}

\renewcommand{\theLemm}{}

\def\a{\alpha}
\def\Aa{{\mathcal A}}
\def\Ss{{\mathcal S}}
\def\b{\beta}
\def\B{{\bf B}} 
\def\C{{\mathcal{C}}} 
\def\CC{{\mathbb{C}}} 
\def\SS{{\mathcal{S}}} 
\def\EE{{\mathbb{E}}} 
\def\Da{{\rm D}_\a}
\def\Dea{\Delta_\a}
\def\Dq{{\rm D}_q}
\def\esp{{\mathbb{E}}} 
\def\elaw{\stackrel{d}{=}}
\def\eps{\varepsilon}
\def\G{{\bf \Gamma}} 
\def\gam{\gamma} 
\def\dt{\delta} 
\def\Ga{\Gamma} 
\def\Fa{F_\a}
\def\Ka{K_\a}
\def\Xa{X_\a} 
\def\ca{c_\a}
\def\fa{f_\a} 
\def\ga{g_\a} 
\def\sa{s_\a} 
\def\WW{{\mathcal{W}}} 
\def\hh{{\hat h}}
\def\hT{{\hat T}}
\def\hX{{\hat X}}
\def\ii{{\rm i}}
\def\K{{\bf K}} 
\def\L{{\mathcal{L}}} 
\def\lb{\lambda}
\def\lacc{\left\{}
\def\lcr{\left[}
\def\lpa{\left(}
\def\lva{\left|}
\def\NN{{\mathbb{N}}} 
\def\pa{p_\a}
\def\Pa{P_\a}
\def\Qa{Q_\a}
\def\Ra{R_\a}
\def\pb{{\mathbb{P}}}
\def\bQa{{\bar Q}_\a}
\def\Qqa{{\widehat \Qa}}
\def\R{{\mathcal{R}}}
\def\rl{{\mathbb{R}}}
\def\racc{\right\}}
\def\rcr{\right]}
\def\rpa{\right)}
\def\rva{\right|}
\def\X{{\bf X}} 
\def\TT{{\rm T}} 
\def\U{{\bf U}} 
\def\Un{{\bf 1}}
\def\Z{{\bf Z}} 

\newcommand{\fin}{\vspace{-0.4cm}
                  \begin{flushright}
                  \mbox{$\Box$}
                  \end{flushright}
                  \noindent}

\title[Total positivity in stable semigroups]{Total positivity in stable semigroups}

\author[Thomas Simon]{Thomas Simon}

\address{Laboratoire Paul Painlev\'e, Universit\'e Lille 1, Cit\'e Scientifique, F-59655 Villeneuve d'Ascq Cedex. Laboratoire de Physique Th\'eorique et Mod\`eles Statistiques, Universit\'e  Paris-Sud, F-91405 Orsay.
{\em Email}: {\tt simon@math.univ-lille1.fr}}



\begin{abstract} We characterize the total positivity in space-time of strictly stable semigroups on $\rl.$ In the positive case, this solves a problem which had been raised by Karlin. In the drifted Cauchy case, this concludes a study which we had initiated in a previous paper. The case of the isotropic stable semigroup on $\rl^d$ is also investigated. We apply these results to the bell-shape and monotone likelihood properties of certain stable densities.
\end{abstract}

\maketitle

\section{Introduction and statement of the results}

The strictly stable semigroup on $\rl$ is a convolution semigroup whose transition density $p_{\a,\rho}(t,x)$ has Fourier transform
\begin{equation}
\label{FT}
\int_\rl p_{\a,\rho}(t,x) e^{\ii\lb x}\, dx\; =\; e^{-t(\ii \lambda)^\a e^{-\ii\pi\a\rho\, {\rm sgn}(\lambda)}},\qquad  t> 0,\, \lb \in\rl.
\end{equation}
The convolution property is meant as $p_{\a,\rho}(t,.)\,\ast\, p_{\a,\rho}(s,.)\, =\, p_{\a,\rho}(t+s, .)$ for all $t,s > 0.$ The parameter $\a\in (0,2]$ is a self-similarity parameter:
\begin{equation}
\label{SSIM}
p_{\a,\rho}(t, x)\; =\; t^{-1/\a} p_{\a,\rho} (1,xt^{-1/\a}), \qquad t> 0,\, x\in\rl.
\end{equation}
Throughout, we will set $f_{\a,\rho} (x) = p_{\a, \rho}(1,x).$ The parameter $\rho\in [0,1]$ is a positivity parameter:
$$\int_0^\infty p_{\a,\rho}(t, x)\, dx\; =\; \rho, \qquad t >0.$$
One has the necessary restrictions $\rho\in[1-1/\a, 1/\a]$ if $\a > 1$ and $\rho \in (0,1)$ if $\a = 1.$ If $\rho=1/2,$ the right-hand side of (\ref{FT}) reads $e^{-t\vert \lambda\vert^\a},$ which is the Fourier transform of the traditional symmetric stable semigroup. If $\a=2,$ then $\rho =1/2$ and 
\begin{equation}
\label{Gauss}
p^{}_{2,1/2}(t,x)\; =\; \frac{e^{-\frac{x^2}{4t}}}{2\sqrt{\pi t}}
\end{equation}
is the standard Gaussian semigroup. It is known that the strictly stable semigroup is up to multiplicative normalization the only convolution semigroup fulfilling a self-similarity property of the type (\ref{SSIM}). We refer e.g. to \cite{ST}, Chapter 8 in \cite{B} and Chapter 3 in \cite{S} for more detail on this classical semigroup and the related stable L\'evy process.

Of particular interest is the positive stable semi-group $p_{\a, 1}(t,x) = p_\a(t,x),$ which is defined for $\a\in(0,1), x >0$ and has a characterization by the Laplace transform:
\begin{equation}
\label{PS}
\int_0^\infty  p_\a(t,x) e^{-\lb x}\, dx\; =\; e^{-t\lambda^\a},\qquad    t,\lb > 0.
\end{equation}
The importance of the positive stable semigroup stems from its subordination properties. For instance, the transition densities of the semigroup $\{e^{-t(-\Delta)^\a}, \,t > 0\}$ built on the fractional Laplacian $-(-\Delta)^\a$ on $\rl^d$ are obtained by the formula
\begin{equation}
\label{Lap}
q_{\a,d}(t,z)\; =\;\int_0^\infty e^{-\frac{\vert\vert z\vert\vert^2}{4x}} \, p_\a(t,x)\, \frac{dx}{(4\pi x)^{d/2}},\qquad t > 0, \, z\in\rl^d.
\end{equation}
A more general result due to Zolotarev shows that every transition density $p_{\a, \rho}(t,x)$ on the line can be recovered by a multiplicative convolution of two positive stable densities - see Formula (3.3.16) in \cite{Z} which will be recalled in (\ref{Facz}) below. Throughout, we will set $f_\a(x) = p_\a(1,x)$ and $\Z_\a$ for the positive random variable with density $f_\a.$ 

Recall that except for $\{\a =1, \rho\in(0,1)\}, \{\a =2, \rho =1/2\}$ and $\{\a=1/2, \rho =1\},$ the densities $p_{\a,\rho}(t,x)$ are not expressed in closed form. Instead, one has convergent series representation given e.g. by (14.30) and (14.31) in \cite{S}. The absence of closed expression makes these densities difficult to study from a classical analytical viewpoint. For instance, it is not even clear from the series representations that the densities are positive.

In this paper, we are interested in the following property extending positivity. Let $I$ be some real interval and $K$ some real kernel defined on $I\times I$. The kernel $K$ is called totally positive of order $n$ (${\rm TP}_n$) if 
\begin{equation}
\label{DET}
\det \lcr K(x_i, y_j)\rcr_{1\le i,j\le m}\; \ge \; 0
\end{equation}
for every $m\in \{1,\ldots, n\},\, x_1< \ldots < x_m$ and $y_1< \ldots < y_m.$ One says that $K$ is ${\rm TP}_\infty$ if these inequalities hold for all $n$. We refer to \cite{K} for the classic account on this field and its numerous connections with analysis. Let us also mention the recent monograph \cite{P} for a more linear algebraic point of view and updated references. General results on the total positivity in space-time of Markovian kernels are displayed in \cite{K} pp. 43-45. The present paper aims at characterizing this property for the kernels $p_{\a,\rho}(t,x)$ with $t,x >0.$ We begin with the positive case.

\begin{Theo} For every $n\ge 2,$ one has 
$$p_\a(t,x)\;\,\mbox{{\em is}}\;\,{\rm TP}_n\,\Leftrightarrow\,\alpha\in \{1/2, 1/3, \ldots, 1/n\}\;\, {\rm or}\;\, \alpha < 1/n.$$ 
 \end{Theo}

This result solves a problem which has some history.  The fact that $p_\a(t,x)$ is ${\rm TP}_\infty$ if and only if $\alpha$ is the reciprocal of an integer is shown in Paragraph 3.4.(vi) pp.121-122 and 7.12.E p.390 of \cite{K}. In the remaining cases, Karlin asks - see again \cite{K} p.390 - whether $p_\a(t,x)$ is  totally positive of some finite order. The fact that $p_\a(t,x)$ is ${\rm TP}_2$ if and only if $\alpha\le 1/2$ is proved in \cite{ThS1}, in connection with the multiplicative strong unimodality of $\Z_\a.$ Finally, the fact that $p_\a(t,x)\in{\rm TP}_n$ implies $\alpha\in \{1/2, 1/3, \ldots, 1/n\}$ or $\alpha < 1/n$ is obtained in the Corollary of \cite{ThS3}. 

From (\ref{SSIM}) and Theorem 1.2.1 in \cite{K}, the ${\rm TP}_n$ character of $p_\a(t,x)$  amounts to that of the convolution kernel $f_\a(xy^{-1}).$ As a rule, the ${\rm TP}_2$ resp. ${\rm TP}_\infty$ property of a convolution kernel is easier to study, because of its characterization through log-concavity - see Theorem 4.1.8 in \cite{K} - resp. Hadamard factorization - see Theorem 7.3.2 in \cite{K}. On the other hand, there is no such handy criterion for testing ${\rm TP}_n$ with $n\ge 3.$ Our method relies on  an infinite multiplicative factorization of the random variable $\Z_\a^{-\a}$ in terms of renormalized Beta random variables. We prove a local limit theorem, which reduces the problem to the total positivity of the fractional integration kernel
$$(x-y)_+^{\frac{1}{\a} -2}.$$
We then apply a result of Karlin on generalized Bessel functions which entails, rather unexpectedly, the right order of total positivity for our kernel $p_\a(t,x).$  

Let us now consider the kernel
$$\Ka(t,x)\; =\; \frac{1}{t^2 + 2 \cos(\pi\a)tx +x^2}\; =\; \frac{\pi}{t\sin(\pi\a)}\,\times\, p_{1,\a}(t,x)$$
over $(0,+\infty)\times(0, +\infty),$ with $\a\in (0,1).$ This can be viewed as a generalization of the classical Cauchy kernel $K_{1/2},$ whose total positivity of infinite order serves as a basic example in \cite{K} p. 149 and is a consequence of Cauchy's double alternant formula - see also Example 4.3 in \cite{P}. As a consequence of the above result, we prove the following non-trivial characterization.

\begin{Corr} For every $n\ge 2,$ one has 
$$K_\a(t,x)\;\,\mbox{{\em is}}\;\,{\rm TP}_n\,\Leftrightarrow\,\alpha\in \{1/2, 1/3, \ldots, 1/n\}\;\, {\rm or}\;\, \alpha < 1/n.$$ 
\end{Corr}

The characterization of the ${\rm TP}_\infty$ property of $\Ka$ through $\alpha\in \{1/2, 1/3, \ldots, 1/n, \ldots\}$ is obtained in our previous paper \cite{ThS3}, thanks to the Hadamard factorization. The fact that $\Ka(t,x)\in{\rm TP}_n$ implies $\alpha\in \{1/2, 1/3, \ldots, 1/n\}$ or $\alpha < 1/n$ is also proven in \cite{ThS3} with the help of an evaluation of the derivative determinant
$$\det \lcr \frac{\partial^{i+j-2}\Ka (t,x)}{\partial t^{i-1}\partial x^{j-1}}\rcr$$
when $t=1$ and $x\to 0,$ by Chebyshev polynomials of the second kind. The reverse inclusion is thoroughly investigated in \cite{ThS3}, viewing the derivative determinant as a generating function of alternating sign matrices. In particular, it is shown that the reverse inclusion holds under a plausible hypothesis on this generating function - see Conjecture 4.1 in \cite{ThS3}. However, this hypothesis is difficult to check. It is surprising to the author that the analysis of the non-explicit kernel $p_\a(t,x)$ turns out to be simpler than that of the explicit kernel $\Ka(t,x).$\\

We now state our main result. Introduce the non-negative parameters 
$$\gamma \;=\; \frac{1}{\rho}\, -\, 1\qquad\quad\mbox{and}\qquad\quad \delta\; =\; \frac{1}{\rho\a}\, -\, 1,$$
and observe that the aforementioned necessary restrictions on $(\a, \rho)$ entail $\gamma\le \delta +1.$ 

\begin{Theo} Assume $\rho \in (0,1)$ and $\a < 2.$ For every $n\ge 2,$ one has 
$$p_{\a,\rho}(t,x)\;\,\mbox{{\em is}}\;\,{\rm TP}_n\,\Leftrightarrow\,(\gam,\dt)\in \NN^2\;\; {\rm or}\;\; \inf(\gam, \dt) \ge n-1.$$ 
\end{Theo}

Notice that the assumption made on $(\a,\rho)$ is no restriction. Indeed, when $\rho=0$ the kernel is zero, whereas the case $\rho =1$ is the matter of Theorem 1 and the case $\a =2$ yields the Gaussian kernel (\ref{Gauss}) which is known to be ${\rm TP}_\infty$ by Formula 3 (1.5) p.100 in \cite{K}. Observe finally that for $\a =1$ the result is already contained in Corollary 1, since $\gam=\dt.$  The following figure illustrates our results, the bold points meaning ${\rm TP}_\infty.$  

\bigskip

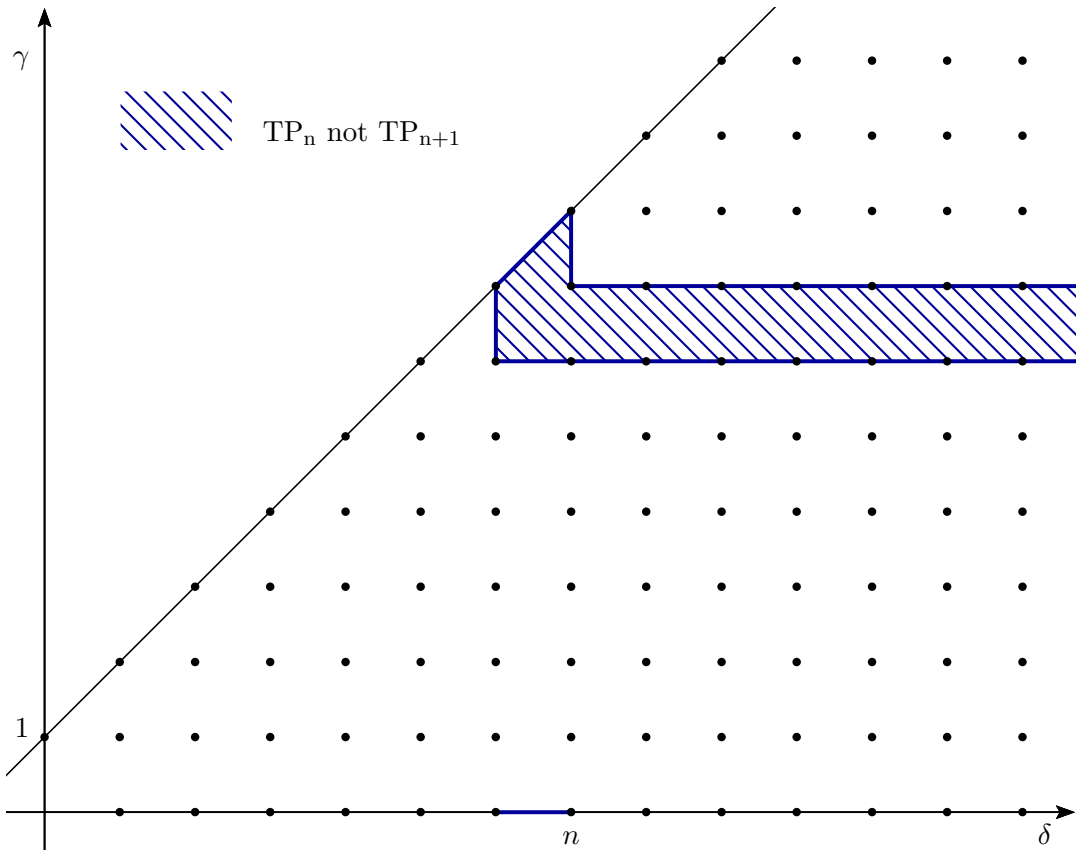
\begin{center}
\definecolor{couleur}{rgb}{0.0,0,0.6}
\definecolor{gris}{rgb}{0.7,0.7,0.7}
\definecolor{couleur2}{rgb}{0,0,0.6}
\centering
\psset{unit=1cm,algebraic=true,linewidth=0.6pt,arrowsize=4pt 2,arrowinset=0.2,subgriddiv=0}
\begin{pspicture*}(-0.5,-0.5)(13.7,10.7)
\psaxes[Dx=30,Dy=15,ticksize=-2pt 2pt]{->}(0,0)(-80,-0.9)(13.7,10.7)
\psline(-5,-4)(20,21)
\pspolygon[hatchcolor=couleur,linewidth=1.5pt,linecolor=couleur,
fillstyle=hlines,hatchangle=-45,hatchsep=0.15](6,7)(7,8)(7,7)(100,7)(100,6)(6,6)
\psline[linecolor=couleur,linewidth=1.5pt](6,0)(7,0)
\newcounter{xpoints}
\newcounter{ypoints}
\newcounter{ymax}
\forloop{xpoints}{1}{\value{xpoints}<20}{
	\setcounter{ymax}{\value{xpoints}+2}
	\forloop{ypoints}{0}{\value{ypoints}<\value{ymax}}{
		\psdots[linecolor=black
		](\arabic{xpoints},\arabic{ypoints})
	}
}

\psdot[linecolor=black](0,1)
\rput[b](-0.3,1){$1$}
\rput[b](7,-0.42){$n$}
\rput[b](13.3,-0.42){$\delta$}
\rput[r](-0.2,10){$\gamma$}
\pspolygon[hatchcolor=couleur,linecolor=white,
fillstyle=hlines,hatchangle=-45,hatchsep=0.15](1,9.6)(2.5,9.6)(2.5,8.8)(1,8.8)
\normalsize
\rput[l](2.9,9.02){$\rm TP_n~not~ TP_{n+1}$} 
\end{pspicture*}
\end{center}

\bigskip

We finally turn our attention to the kernel 
$$\rho_{\a,d}(t,r)\; =\;\int_0^\infty e^{-\frac{r^2}{4x}} \, p_\a(t,x)\, \frac{dx}{(4\pi x)^{d/2}},\qquad t, r > 0,$$
which is the radial part of the transition density $q_{\a,d}(t,z)$ on $\rl^d.$ The following can be viewed as a multidimensional generalization of Theorem 2, in the case $\rho = 1/2.$ 

\begin{Corr}
For every $n\ge 2,$ one has 
$$\rho_{\a,d}(t,r)\;\,\mbox{{\em is}}\;\,{\rm TP}_n\,\Leftrightarrow\,\alpha\in \{1/2, 1/3, \ldots, 1/n\}\;\, {\rm or}\;\, \alpha < 1/n.$$ 
\end{Corr}

The if part of the proof of Theorem 2 resp. Corollary 2 is a rather straightforward consequence of Theorem 1 and Zolotarev's aforementioned factorization resp. the chi-square factorization which is obtained from (\ref{Lap}). The only if part is more delicate and relies on an asymptotic analysis of the derivative determinant, in the same spirit as our previous papers \cite{ThS2, ThS3}. These proofs are given in the next section after the proofs of Theorem 1 and Corollary 1. In the third section, we apply some of these results to two visual features of real stable densities: the bell-shape and the monotone likelihood property.

\section{Proofs}

\subsection{Proof of Theorem 1} 

From the discussion in the introduction, we need to show the implication $\a < 1/n \Rightarrow f_\a(x y^{-1})$ is ${\rm TP}_n$ for every $n\ge 2.$ Our argument relies on an analysis of the proximity between the random variable $\X_\a = \Z_\a^{-\a}$ and the independent product
$$\X_{\a, n}\; =\; e^{\gamma(\a-1)}\,\times\,\prod_{k=0}^n e^{\psi(1+\frac{k}{\a}) - \psi(\frac{1+k}{\a})}\, \B_{1+\frac{k}{\a},\frac{1}{\a}-1}$$
where $\gamma$ is Euler's constant, $\psi$ is the digamma function and $\B_{a,b}$ is the usual Beta random variable with density 
$$f_{\B_{a,b}}^{}(x)\; =\;\frac{\Ga(a+b)}{\Ga(a)\Ga(b)} \, x^{a-1}(1-x)^{b-1} \, \Un_{(0,1)}(x).$$
Recall the expression of the fractional moments of $\B_{a,b}$:
\begin{equation}
\label{Frac}
\esp[\B_{a,b}^s]\; =\; \frac{\Ga(a+s)\Ga(a+b)}{\Ga(a)\Ga(a+b+s)}, \qquad s > -a.
\end{equation}
We first show the following local limit theorem.
\begin{Lemm} Let $g_\a$ and $g_{\a,n}$ be the respective densities of $\X_\a$ and $\X_{\a, n}.$ One has
$$\sup_{x > 0}\vert g_\a(x) - g_{\a, n} (x)\vert\; \to\; 0 \qquad \mbox{{\em as} $n \to +\infty.$}$$
\end{Lemm}

\proof Setting $h_\a(x) = e^x g_\a(e^x)$ and $h_{\a,n}(x) = e^x g_{\a,n}(e^x)$ for the respective densities of $\log \X_\a$ and $\log \X_{\a, n},$ we have to prove that
\begin{equation}
\label{Exp}
\sup_{x \in \rl}\vert h_\a(x) - h_{\a, n} (x)\vert\; \to\; 0 \qquad \mbox{as $n \to +\infty.$}
\end{equation}
Let ${\hh_\a}$ and ${\hh_{\a,n}}$ be the Fourier transforms of $h_\a$ and $h_{\a,n}.$ On the one hand, we have e.g. from Formula 3.(4.12) p.121 in \cite{K}
$${\hh_\a}(s)\; =\;\esp[\Z_\a^{-\ii\a s}]\; =\; \frac{\Ga(1+\ii s)}{\Ga(1+\ii\a s)}$$
which is an integrable function over $\rl$ by Formula 1.18(6)  in \cite{EMOT}, since $\a\in(0,1).$ On the other hand, supposing $n\ge k_\a =\inf\{k\ge 1, k(1-\a) > \a\},$ it follows from (\ref{Frac}) that
$$\vert {\hh_{\a,n}}(s)\vert\; =\; \prod_{k=0}^n \lva \frac{\Ga(1+\frac{k}{\a} +\ii s)\Ga(\frac{1+k}{\a})}{\Ga(\frac{1+k}{\a}+\ii s)\Ga(1+\frac{k}{\a})}\rva\;\le\; \prod_{k=0}^{k_\a} \lva \frac{\Ga(1+\frac{k}{\a} +\ii s)\Ga(\frac{1+k}{\a})}{\Ga(\frac{1+k}{\a}+\ii s)\Ga(1+\frac{k}{\a})}\rva\; =\; O(\vert s\vert^{-(k_\a+1)(\frac{1}{\a} -1)}),$$
where the inequality comes from the fact that all multiplicands are characteristic functions, and the estimate at infinity follows at once from Stirling's formula - see e.g. Formula 1.18(4)  in \cite{EMOT}. This entails that ${\hh_{\a,n}}$ is integrable over $\rl$ for all $n\ge k_\a,$ and we can apply the Fourier inversion formula to get the uniform bound
$$\sup_{x \in \rl}\vert h_\a(x) - h_{\a, n} (x)\vert\; \le \; \frac{1}{2\pi}\int_\rl \vert  {\hh_\a}(s) - {\hh_{\a,n}}(s)\vert \,ds.$$
It is clear from the above discussion that
$$\sup_{n\ge k_\a}\lpa\int_{[-K,K]^c} \vert  {\hh_\a}(s) - {\hh_{\a,n}}(s)\vert\,ds\rpa\; \le \; \int_{[-K,K]^c}(\vert  {\hh_\a}(s)\vert + \vert {\hh_{\a,k_\a}}(s)\vert)\,ds\; \to \; 0\qquad\mbox{as $K\to +\infty$.}$$
By dominated convergence, it is hence sufficient to show that $\hh_\a \to \hh_{\a, n}$ pointwise as $n\to +\infty.$ Proceeding as in Lemma 2 of \cite{BS}, for every $s\in\rl$ we observe that
\begin{eqnarray*}
\hh_\a(s) & =&  \exp\lcr \ii\gamma(\a-1)s\; +\; \int_{-\infty}^{0} (e^{\ii\a sx} - 1 -\ii\a sx)\, \frac{e^{-\a\vert x\vert} (1-e^{-(1-\a)\vert x\vert})}{\vert x\vert (1-e^{-\vert x\vert}) (1-e^{-\a\vert x\vert})}\, dx\rcr \\
& = & \exp\lcr \ii s\gamma(\a-1)\; +\;\sum_{k=0}^\infty \lpa\int_{-\infty}^{0} (e^{\ii sx} - 1 -\ii sx)\, \frac{e^{-(1+\frac{k}{\a})\vert x\vert} (1-e^{-(\frac{1}{\a} -1)\vert x\vert})}{\vert x\vert (1-e^{-\vert x\vert})}\, dx\rpa\rcr\\
& = & \lim_{n\to +\infty}\lpa \exp\lcr \ii s\gamma(\a-1)\; +\;\sum_{k=0}^n\lpa \int_{-\infty}^{0} (e^{\ii sx} - 1 -\ii sx)\, \frac{e^{-(1+\frac{k}{\a})\vert x\vert} (1-e^{-(\frac{1}{\a} -1)\vert x\vert})}{\vert x\vert (1-e^{-\vert x\vert})}\, dx\rpa\rcr\rpa\\
& = & \lim_{n\to +\infty}\lpa e^{\ii s\gamma(\a-1)}\,\times\, \prod_{k=0}^n e^{\ii s(\psi(1+\frac{k}{\a}) - \psi(\frac{1+k}{\a}))}\, \frac{\Ga(1+\frac{k}{\a} +\ii s)\Ga(\frac{1+k}{\a})}{\Ga(\frac{1+k}{\a}+\ii s)\Ga(1+\frac{k}{\a})}\rpa \; = \; \lim_{n\to +\infty} \hh_{\a,n}(s).
\end{eqnarray*}
This completes the proof.

\endproof

We can now finish the proof of Theorem 1. Again by Theorem 1.2.1 in \cite{K}, we must show that the kernel
$$g_\a(xy^{-1})\; =\; \frac{y^{1+1/\a}}{\a x^{1+1/\a}}\, f_\a(y^{1/\a}x^{-1/\a})$$
is ${\rm TP}_n$ when $\a < 1/n.$ By Lemma 1 and a straightforward continuity argument, it is enough to prove the same property for $g_{\a,k}(xy^{-1})$ for every $k\ge 0.$ This kernel is of the multiplicative convolution type and applying repeatedly Lemma 3.1.1 in \cite{K}, we see that it is sufficient to prove that
$$f_{\B_{1+\frac{k}{\a},\frac{1}{\a}-1}}^{}(xy^{-1})\; =\; \frac{\Ga(\frac{1+k}{\a})\,x^{\frac{k}{\a}}\,y^{2-\frac{1+k}{\a}}}{\Ga(1+\frac{k}{\a})\Ga(\frac{1}{\a}-1)} \; (y-x)_+^{\frac{1}{\a}-2}$$ 
is ${\rm TP}_n$ for $\a < 1/n$ and every $k\ge 0,$ with the notation $z_+ = z \Un_{\{z\ge 0\}}.$ The latter property is a direct consequence of Theorem 3.2.1 in \cite{K}, making $c_0 = 1$ and $c_n =0, n\neq 0$ therein.

\qed

\begin{Rq} {\em A similar characterization of the ${\rm TP}_n$ property can be observed for the fractional integration kernel $I_\b(x,y) = (y-x)_+^{\b-1}$ with $\b > 0.$ One has
$$I_\b(x,y)\,\in{\rm TP}_{n}\,\Leftrightarrow\,\b\in \NN\;\; {\rm or}\;\; \b > n-1.$$ 
The if part follows from Theorem 10.1.1 and Theorem 3.2.1 in \cite{K}, whereas the only if part can be deduced from Part (ii) of the main result in \cite{ThS2}, making $d\le 0$ therein. We omit the details. }

\end{Rq}

\subsection{Proof of Corollary 1} As above, we only need to show the reverse inclusion that is $\K_\a\in{\rm TP}_n$ if $\a < 1/n,$ for every $n\ge 2.$ To do so, we appeal to a standard computation which was already used in \cite{ThS1} - see (3.1) and the references therein, and which reads
$$\Ka(e^{\a x}, e^{\a y})\; =\; \frac{\pi e^{(1-\a)(x+y)}}{\sin(\pi\a)} \int_\rl  \fa(e^{x-u})\fa(e^{y-u})\,e^{-2u}\, du$$
for every $x,y\in \rl.$ It is clear from Theorem 1 and Lemma 3.1.1 in \cite{K} that the integral on the right-hand side defines a ${\rm TP}_n$ kernel on $\rl\times\rl$ as soon as $\a < 1/n.$ This completes the proof, again by Theorem 1.2.1 in \cite{K}.

\qed

\begin{Rq} {\em (a) The discrete set characterizing the total positivity of $\Ka$ is in one-to-one correspondence with the set $\{\cos(\pi/(n+1)),\, n\ge 1\}$ of the largest roots of the Chebyshev polynomials of the second kind 
$$U_n(\cos\theta)\; =\;\frac{\sin(n+1)\theta}{\sin\theta}\cdot$$ 
This is not a surprise, recalling the formula  
$$\sum_{n\ge 0} (-1)^n z^n \,U_n(\cos\pi\a) \; =\;\Ka(1,z)$$
for the generating function. On the other hand, plugging this generating function into determinants yields a complicated alternate series which seems unappropriate to tackle the total positivity problem - see Remark 3.2 in \cite{ThS3}. Notice that zeroes of Chebyshev polynomials appear in connection with the total positivity of infinite Jacobi matrices - see \cite{K} pp.115-117. Observe also the connection with Katkova-Vishnyakova's general criterion on $2\times 2$ minors - see Theorem 2.16 in \cite{P}, which is however too stringent in our framework.\\

(b) The main result of \cite{ThS3} entails the further equivalence
$$K_\a(t,x)\;{\rm is}\;{\rm SR}_n\,\Leftrightarrow\,K_\a(t,x)\;{\rm is}\;{\rm TP}_n\,\Leftrightarrow\,\alpha\in \{1/2, 1/3, \ldots, 1/n\}\;\, {\rm or}\;\, \alpha < 1/n,$$
where the left-hand side means sign-regularity of order $n$ - see \cite{K} p.12. It is conjectured in \cite{ThS3} that the ${\rm TP}_n$ property of $K_\a(t,x)$ is actually a ${\rm STP}_n$ property, that is all the determinants of (\ref{DET}) should be positive. This is true for $\alpha\in \{1/2, 1/3, \ldots, 1/n, \ldots\}$ by Proposition 2.1 in \cite{ThS3}. In the remaining case $\a < 1/n,$ it does not seem that the factorization methods of the present paper can help to show this stronger property, because the fractional integration kernels $I_\beta$ are never ${\rm STP}_n$. See Section 4 in \cite{ThS3} for  a partial result and a related conjecture on the generating function of ASM matrices with fixed number of negative entries.}
\end{Rq}

\subsection{Proof of Theorem 2} Again by (\ref{SSIM}), the ${\rm TP}_n$ property of $p_{\a,\rho}$ amounts to that of the convolution kernel $f_{\a,\rho}(xy^{-1}),$ where $f_{\a,\rho}$ is the density of the positive part $\X^+_{\a,\rho}$ of the stable random variable with parameters $(\a,\rho).$ By Formula (3.3.16) in \cite{Z}, one has the independent factorization
\begin{equation}
\label{Facz}
\X^+_{\a,\rho}\;\elaw\;\lpa\frac{\Z_{\a\rho}}{\Z_\rho}\rpa^\rho,
\end{equation}
where we have made the convention $\Z_1 \equiv \Un.$ If $(\gam,\dt)\in \NN^2$ or $\inf(\gam, \dt) \ge n-1,$ we see from Theorem 1 that the kernels $f_{\a}(xy^{-1})$ and $f_{\a\rho}(xy^{-1})$ are ${\rm TP}_n.$ By (\ref{Facz}), Lemma 3.1.1 in \cite{K} and the multiplicative convolution formula, we deduce that the same holds for $f_{\a,\rho}(xy^{-1}),$ which concludes the if part of Theorem 2.\\

We now turn to the only if part. Recall that the case $\a =1$ is already worked out in Corollary 1. Besides, the identity in law
$$\X^+_{\a,\rho}\;\elaw\;(\X^+_{\frac{1}{\a},\a\rho})^{-\frac{1}{\a}}$$
for every $\a\in (1,2)$ and $\rho\in[1-1/\a, 1/\a],$ which is known as Zolotarev's duality formula and is a direct consequence of (\ref{Facz}), reduces the study to the case $\a < 1,$  which we henceforth assume. We will reason as in Section 3 of \cite{ThS3}, using the smoothness of the function $f_{\a, \rho}$ and an analysis on the derivative determinant. Supposing $f_{\a,\rho}(xy^{-1})\in{\rm TP}_n$ entails the everywhere non-negativity of 
$$(x,y)\;\mapsto\; (-1)^{\frac{k(k-1)}{2}}\,{\det}_k \lcr \frac{\partial^{i+j-2}K_{\a,\rho} (x,y)}{\partial x^{i-1}\partial y^{j-1}}\rcr$$
for every $k=2,\ldots, n,$ where we have set $K_{\a,\rho} (x,y) = f_{\a,\rho}(xy)$ and ${\det}_k$ means the determinant of a $k\times k$ matrix. Elementary computations lead then to 
\begin{equation}
\label{DetP}
\Delta^k_{\a,\rho}(z)\; =\; (-1)^{\frac{k(k-1)}{2}}\,{\det}_k \lcr \lpa z^{j-1} f_{\a,\rho}^{(j-1)}(z)\rpa^{(i-1)} \rcr\; \ge \; 0
\end{equation}
for every $k=2,\ldots, n$ and $z \ge0,$ where $f^{(q)}$ means the $q-$th derivative of a given function $f.$ \\

Suppose first $\gam < n-1$ and $\gam\not\in\NN.$ Formula (14.33) in \cite{S} shows that
$$ f_{\a,\rho}^{(j-1)}(0)\; =\; \frac{(-1)^{j-1} \Ga(1+\frac{j}{\a})}{\pi j}\, \times\, \sin(\pi j\rho)\;\neq\; 0$$
for every $j = 1, \ldots, n.$ Simplifying the determinant, this implies that
$$\Delta^k_{\a,\rho}(0)\; =\; \lpa \frac{\prod_{j=1}^k \Ga(1+\frac{j}{\a})}{\pi^k k!} \rpa\,\times\, \lpa\prod_{j=1}^k \,\sin(\pi j\rho)\rpa$$
and is necessarily negative for some $k\in \{2,\ldots, n\},$ a contradiction.\\

Suppose next $\dt < n-1$ and $\dt\not\in\NN.$ We will let $z\to +\infty$ and obtain a similar contradiction, but the argument is more subtle. Fix $k\in \{2,\ldots, n\}$ and set $g_j(z) = (-1)^{j-1} z^{j-1} f_{\a,\rho}^{(j-1)}(z)$ for every $j=1,\ldots, k.$ With this notation, we have the Wronskian representation
$$\Delta^k_{\a,\rho}\; =\; \WW(g_1, \ldots, g_k).$$
On the other hand, Formula (14.31) in \cite{S} yields the convergent series representation 
$$g_j^{(i-1)}(z)\; =\; \frac{(-1)^{i-1}}{\pi} \sum_{q=1}^\infty (-1)^{q-1} \frac{\Ga(j+q\a)}{q!}\lpa\prod_{r=1}^{i-1}\,(r+q\a)\rpa\, \times\, \sin(\pi q\rho\a)\, z^{-q\a -i}, \qquad z >0.$$
We see that the terms of fixed order $q$ in the expansion of the lines
$$(g_1^{(i-1)}(z)\ldots g_k^{(i-1)}(z))$$
are linearly dependent, and this entails by multilinear expansion that the leading term in the expansion of $\Delta^k_{\a,\rho}(z)$ as $z \to +\infty$ is obtained in taking terms of different order from $1$ to $k$ in each line: we get
$$\Delta^k_{\a,\rho}(z)\; \sim\;\frac{1}{\pi^k}\lpa\sum_{\sigma\in\SS_k}{\det}_k \lcr\frac{\Ga(j+\sigma(i)\a)}{\sigma(i)!}\lpa\prod_{r=1}^{i-1}\,(r+\sigma(i)\a)\rpa\, \times\, \sin(\pi \sigma(i)\rho\a)\rcr\rpa z^{-\frac{k(k+1)(\a+1)}{2}}$$
and see that the leading coefficient is, up to a positive constant, given by the product 
$$\lpa\sum_{\sigma\in\SS_k}{\det}_k \lcr\Ga(j+\sigma(i)\a)\lpa\prod_{r=1}^{i-1}\,(r+\sigma(i)\a)\rpa\rcr\rpa\,\times\,\lpa \prod_{j=1}^k\, \sin(\pi j\rho\a)\rpa.$$
Hence, to obtain the same contradiction as above, it suffices to show that
$$\sum_{\sigma\in\SS_k}{\det}_k \lcr\Ga(j+\sigma(i)\a)\lpa\prod_{r=1}^{i-1}\,(r+\sigma(i)\a)\rpa\rcr\; > \; 0$$
for every $k\ge 2.$ Making the necessary simplifications, we transform the latter expression into
\begin{eqnarray*}
\sum_{\sigma\in\SS_k} (-1)^{\eps(\sigma)} {\det}_k \lcr\Ga(i +\a j)\rcr\times \prod_{i=1}^k\lpa\prod_{r=1}^{i-1}\,(r+\sigma(i)\a)\rpa & = & {\det}_k \lcr\Ga(i +\a j)\rcr \,\times\, {\det}_k \lcr \prod_{r=1}^{i-1}\,(r+\a j) \rcr\\
& = & \lpa\prod_{j=1}^k \Ga(1 +\a j) \rpa\times\, {\det}_k^2 \lcr \prod_{r=1}^{i-1}\,(r+\a j) \rcr\\
& = & \!\a^{k(k-1)} \lpa\prod_{j=1}^k  ((j-1)!)^2 \,\Ga(1 +\a j)\rpa\, > \, 0.
\end{eqnarray*}
This completes the proof.

\qed

\begin{Rq} {\em This results entails $p_{\a,\rho}(t,x)\,\mbox{is}\,{\rm TP}_\infty\,\Leftrightarrow\,(\gam,\dt)\in \NN^2,$ which can be shown directly from the fractional moment evaluation
$$\EE[(\X^+_{\a,\rho})^s]\; =\; \frac{\sin (\pi\rho s)}{\rho \sin (\pi s)}\,\times\, \frac{\Ga(1-\frac{s}{\a})}{\Ga(1-s)}$$
given in Formula (2.6.20) of \cite{Z}, and the Hadamard factorization criterion given in  Theorem 7.3.2 of \cite{K}.}
\end{Rq}

\subsection{Proof of Corollary 2} By a change of variable, the ${\rm TP}_n$ property of $\rho_{\a, d}(t,r)$ amounts to that of the kernel
$$\frac{1}{\Ga(\frac{d}{2})} \int_0^\infty \lpa \frac{r}{tx}\rpa^{\frac{d}{2} -1} e^{-\frac{r}{tx}} \,\fa(x)\, \frac{dx}{x}\; =\; f_{\a, d} (r t^{-1}),$$
where $\G_a$ is the Gamma random variable with density 
$$\frac{1}{\Ga(a)}\, x^{a-1} e^{-x} \Un_{(0,+\infty)}(x)$$
and $f_{\a,d}$ is the density of the independent product $\G_{\frac{d}{2}} \times \Z_\a.$ Since the multiplicative convolution kernel of $\G_{\frac{d}{2}}$ is clearly ${\rm TP}_\infty,$ the if part of Corollary 2 follows again from Lemma 3.1.1 in \cite{K} and Theorem 1.

To show the only if part, we suppose first that $d$ is odd. If $d =1,$ the conclusion follows from Theorem 2. If $d =2p+1\ge 3,$ we use the elementary independent factorization
$$\G_{\frac{1}{2}}\, \times\, \Z_\a\; =\; \lpa\G_{\frac{d}{2}}\, \times\, \Z_\a\rpa\,\times\, \B_{\frac{1}{2}, p}$$
which follows from a fractional moment identification. Suppose $\a > 1/n$ and $\a \not\in\{1/2, \ldots, 1/n\}.$ If $f_{\a, d} (r t^{-1})$ were ${\rm TP}_n,$ then so would be $f_{\a, 1} (r t^{-1})$ by this factorization, Lemma 3.1.1 in \cite{K} and Remark 1. However, this is false and we have a contradiction. Similarly, to show the only if part when $d=2p$ is even, we are reduced to the case $d =2.$

Integrating (14.31) of \cite{S} along the exponential kernel yields the well-known formula
$$f_{\a, 2} (z)\; =\; \sum_{q=1}^\infty (-1)^{q-1} \frac{\Ga(1+q\a)^2}{q!}\, \times\, \sin(\pi q\a)\, z^{-q\a -1}, \qquad z >0,$$
which is meant as a convergent series representation for $\a < 1/2$ and as an asymptotic expansion at infinity for $\a \ge 1/2.$ One can check that this formula is differentiable term by term if $\a \ge 1/2.$  Reasoning as in the proof of Theorem 2 leads to an asymptotic analysis of the Wronskians $\WW(h_1, \ldots, h_k),$
with 
$$h_j^{(i-1)}(z)\; =\; \frac{(-1)^{i-1}}{\pi} \sum_{q=1}^\infty (-1)^{q-1} \frac{\Ga(j+q\a)\Ga(1+q\a)}{q!}\lpa\prod_{r=1}^{i-1}\,(r+q\a)\rpa\, \times\, \sin(\pi q\a)\, z^{-q\a -i},$$
this formula being again meant as an asymptotic expansion for $\a \ge 1/2.$ The analysis is exactly the same as the above and the result follows from the positivity of
$$\sum_{\sigma\in\SS_k}{\det}_k \lcr\Ga(1+\sigma(i)\a)\Ga(j+\sigma(i)\a)\lpa\prod_{r=1}^{i-1}\,(r+\sigma(i)\a)\rpa\rcr\; = \; \a^{k(k-1)} \lpa\prod_{j=1}^k (j-1)! \,\Ga(1 +\a j)\rpa^2$$
for every $k\ge 2.$

\qed

\section{Applications to stable densities}

In this section, we investigate further analytical properties of the real stable densities 
$$f_{\a, \rho}(x)\;=\;p_{\a,\rho} (1,x).$$ 
Let us again refer to \cite{Z} for a thorough presentation of these densities, which are of constant use in probability and statistics. In Zolotarev's opinion, stable densities should also be accorded civil rights in the theory of special functions - see Section 2.10 in \cite{Z}. However, in the absence of closed formul\ae, some of their basic properties are still unknown, or at least unproven.
 
\subsection{Bell-shape} A smooth density function on a real interval whose derivatives vanish at both ends of the interval is said to be {\em bell-shaped} if its $n$-th derivative vanishes exactly $n$ times inside the interval, for every $n\ge 0.$ In the following, we will use the notation ${\rm BS}_n$ if the property holds for every $k$-th derivative with $k\le n.$ For example, a ${\rm BS}_2$ density is strictly positive, strictly unimodal and successively convex, concave and convex like the familiar bell curve. There is a tight connection between bell-shape and total positivity which is exemplified by the notion of ETP kernel - see Paragraph 6.11.C in \cite{K}. In a previous paper \cite{ThS4}, to which we refer for more details, references and open questions around the bell-shape, we have shown that all densities $f_\a$ are in ${\rm BS}_\infty.$

In this paragraph, we apply our previous results to the ${\rm BS}_n$ property of the symmetric density $q_\a(x)= q_{\a,1}(x)=f_{\a, 1/2}(x).$ It is easily seen by Fourier inversion that $$q_\a(x)\; =\; \frac{1}{\pi}\int_0^\infty \cos(tx) \,e^{-t^{2\a}}\, dt $$ 
is real analytic with all derivatives vanishing at infinity. Moreover it is strictly positive on $\rl$ and strictly unimodal by Yamazato's theorem - see Theorem 53.1 in \cite{S} and also \cite{ThS5} for a short proof, hence it is ${\rm BS}_1.$ Let us prove the following reinforcement.

\begin{Prop} For every $n\ge 2,$ the function $q_\a$ is ${\rm BS}_n$ if $\alpha\in \{1/2, 1/3, \ldots, 1/n\}$ or $\alpha < 1/n.$
\end{Prop}

\proof By an induction and a successive application of Rolle's theorem, it is clear that $q_\a^{(n)}$ vanishes and changes its sign at least $n$ times on $\rl,$ since $q_\a^{(n-1)}$ vanishes at infinity. On the other hand, differentiating (\ref{Lap}) shows after a change of variable that
$$q_\a^{(n)}(x)\; =\; \frac{(-1)^n}{\sqrt{2\pi}}\int_0^\infty e^{-x^2 z^{-2}} H_n (xz^{-1})\, p_\a(z^2/4)\, z^{-n} dz$$
for all $x\in \rl, n\ge 0,$ where $H_n$ is the $n-$th Hermite polynomial. In particular, for every $p\ge 0$ one has $q_\a^{(2p+1)}(0)=0.$ Recalling that $q_\a$ is even and putting everything together, we are hence reduced to show that for every $n\ge 2$, the function $q_\a^{(n)}$ vanishes at most $[\frac{n+1}{2}]$ times on $(0,+\infty)$ as soon as $\alpha\in \{1/2, 1/3, \ldots, 1/n\}$ or $\alpha < 1/n.$

Fix $n\ge 2$ and assume $\alpha\in \{1/2, 1/3, \ldots, 1/n\}$ or $\alpha < 1/n.$ A further change of variable implies
\begin{equation}
\label{Ch}
q_\a^{(k)}(x)\; =\; \frac{(-1)^n}{\sqrt{2\pi}}\int_0^\infty e^{-z^2} z^{k-2} H_k (z)\,p_\a(x^2/4z^2)\, dz
\end{equation}
for all $x, k > 0.$ It is known that the function $e^{-z^2} z^{n-2} H_n (z)$ vanishes exactly $[\frac{n+1}{2}]$ times on $(0,+\infty).$ Since $[\frac{n+1}{2}]\le n-1,$ we can deduce from (\ref{Ch}), Theorem 1 and Schoenberg's variation-diminishing property - see Theorem 1.3.1(a) in \cite{K} - that $q_\a^{(n)}$ changes its sign at most $[\frac{n+1}{2}]$ times inside $(0,+\infty),$ the zero terms being discarded. Suppose now that $q_\a^{(n)}$ would vanish on $(0,+\infty)$ without changing its sign. By analyticity, this zero would be isolated and we can again apply Rolle's theorem to deduce that $q_\a^{(n+1)}$ would change its sign at least $1+[\frac{n+2}{2}]$ times. However, since $[\frac{n+2}{2}]\le n-1,$ this is a contradiction because the same argument as above shows that $q_\a^{(n+1)}$ changes its sign at most $[\frac{n+2}{2}]$ times. This completes the proof.

\endproof

\begin{Rq} {\em A direct analysis shows that the Cauchy density
$$q_{1/2,1}(x)\; =\; f_{1,1/2}(x)\; =\; \frac{1}{\pi(1+x^2)}$$
is ${\rm BS}_\infty,$ and by translation the same property holds for all densities $p_{1,\rho}.$ It is very natural to conjecture that all densities $f_{\a,\rho}$ are ${\rm BS}_\infty.$ This is claimed in \cite{Ga}, with unfortunately a serious mistake in the proof - see Section 3.3 in \cite{ThS4} for details. The above proposition is based on the subordination method and can be extended to some densities $f_{\a,\rho}$ with $\a < 1,$ but Theorem 1 clearly shows the limits of this method. To prove the conjecture, an additive factorization like the one used in \cite{ThS4} could be the right approach, eventhough the situation is more complicated on the line than on the half-line.}
 
\end{Rq}

\subsection{Monotone likelihood property} An absolutely continuous positive random variable $X$ is said to have the single intersection property if the density functions of $X$ and $cX$ meet exactly once for every $c\in (0,1)\cup (1,\infty).$ This can be viewed as a refinement as the usual stochastic ordering between $X$ and $cX.$ An absolutely continuous real random variable $Y$ is said to have the double intersection property if the density functions of $Y$ and $cY$ meet exactly twice for every $c\in (0,1)\cup (1,\infty).$ In the symmetric integrable case, this implies a convex ordering between $Y$ and $cY.$ In the following, we will use the same intersection terminology for a given density as for a random variable. We refer to \cite{Sh} for more information on intersection properties for densities and stochastic orderings between random variables.

It was shown in \cite{K} that $f_\a$ has the single intersection property and that $q_\a$ has the double intersection property - see Theorems 4.1 and 4.2 therein. The following proposition, whose proof is independent of the previous results, completes the picture.

\begin{Prop} The density $f_{\a, \rho}$ has the double intersection property as soon as $\rho \in (0,1).$ 
\end{Prop}

\proof With our above notation, it is enough to show that both positive random variables $\X^+_{\a,\rho}$ and $\X^+_{\a,1-\rho}$ have the single intersection property. Indeed, the first one is distributed according as the positive branch of $f_{\a, \rho}$, whereas the second one is distributed according as the negative branch of $f_{\a, \rho}.$ Since both branches are analytic on $(0,+\infty),$ the principle of isolated zeroes and a perusal of the proofs of Lemma 3.1 and Theorem 3.1 imply that we are reduced to prove the unimodality of the random variables $\log\X^+_{\a,\rho}$ and $\log\X^+_{\a,1-\rho}.$ By (\ref{Facz}) and the aforementioned Yamazato theorem, it is enough to show that $\log \Z_\b$ is self-decomposable for every $\b\in (0,1).$ This well-known fact is a consequence of the Fourier transform
$$\EE[e^{\ii s \log \Z_\b}]\;=\;\exp\lcr \ii\gamma(1-1/\b)s\; +\; \int_{-\infty}^{0} (e^{\ii sx} - 1 -\ii  sx)\, \frac{e^{-\b\vert x\vert} (1-e^{-(1-\b)\vert x\vert})}{\vert x\vert (1-e^{-\vert x\vert}) (1-e^{-\b\vert x\vert})}\, dx\rcr $$
which was used in the proof of Theorem 1, the easily established increasing character of
$$x\;\mapsto\; \frac{e^{-\b\vert x\vert} (1-e^{-(1-\b)\vert x\vert})}{(1-e^{-\vert x\vert}) (1-e^{-\b\vert x\vert})}$$
on $(-\infty,0),$ and Corollary 15.11 in \cite{S}.

\endproof

Following Section 5 in \cite{K}, a real density function $f$ is said to have the monotone likelihood ratio (MLR) property in $(0,+\infty)$ if the function $x\mapsto f(x)/f(cx)$ is  monotone on $(0,+\infty)$ for every $c\in (0,1)\cup (1,\infty).$ In the case of strict monotonicity, this property is clearly a reinforcement of the single intersection property. The Karlin-Rubin theorem - see e.g. Theorem 8.4.1 in \cite{M} - establishes a well-known connection between the MLR property and the construction of certain UMP tests. It is observed in Section 5 \cite{K} that none of the functions $f_\a$ and $q_\a$ has the MLR property in general. We can show the following characterization, which somehow contradicts the introduction in \cite{K}.

\begin{Prop} The density $f_{\a, \rho}$ has the {\em MLR} property in $(0,+\infty)$ if and only if one of the three following disjoint conditions is fulfilled: $\{\inf(\gam,\dt) \ge 1\},\, \{\gam = 0, \;\dt \ge 1\},\,\{\gam = 1,\; \dt =0\}.$ 
\end{Prop}

\proof As observed in the beginning of Section 5 in \cite{K}, the MLR property of a density $f$ on $(0,+\infty)$ is equivalent to the ${\rm TP}_2$ character of the kernel $f(xy^{-1})$ on $(0,+\infty)\times(0,+\infty).$ The result is hence a direct consequence of Theorems 1 and 2.

\endproof

In the positive case $\rho =1,$ the MLR property of $\fa$ is hence characterized by $\a\le 1/2,$ which is the main result of \cite{ThS1}. Since $\fa$ is smooth, it is easily proved - see again \cite{ThS1} - that the MLR property of $\fa$ amounts to the non-increasing character of
\begin{equation}
\label{Yam}
x\;\mapsto\; \frac{x\fa'(x)}{\fa(x)}
\end{equation}
on $(0,+\infty).$ It is interesting to ask how this function varies when $\a > 1/2,$ and we believe that it decreases and then increases. The following last proposition provides a partial answer.

\begin{Prop} When $\a > 1/2,$ the function in {\em (\ref{Yam})} decreases on $(0, m_\a],$ where $m_\a$ is the unique positive mode of $f_\a.$
\end{Prop}

\proof We only sketch the proof, which depends on some of the material displayed in Chapter 10 of \cite{S}. Differentiating three times Steutel's integro-differential equation, we have 
$$(x \fa''(x) + \fa'(x)) \fa(x) \, - \, x(\fa'(x))^2\; =\; \frac{\a\fa(x)}{\Ga(1-\a)}\int_0^\infty \fa(x-y) \lpa \frac{\fa'(x)}{\fa(x)} - \frac{\fa'(x-y)}{\fa(x-y)}\rpa \,y^{-\a}\, dy$$
for every $x >0$. See Formula (53.13) in \cite{S}, which is proved there only for a certain class of self-decomposable densities, but it is easily extended to $f_\a$ by approximation. On the other hand, we know that $\fa$ is log-concave on $(0,m_\a]$ by the Yamazato property - see Lemma 53.2 (iii) in \cite{S} and observe that this is again easily extended to $\fa.$ This implies that the right-hand side in the above equation is negative for every $x\in(0,m_\a].$ Hence, the function in (\ref{Yam}) decreases on $(0, m_\a].$

\endproof

\begin{Rq} {\em This proposition implies that the function $\fa(x)/\fa(cx)$ increases on $(0, m_\a/c]$ for any $c > 1$ and decreases on $(0,m_\a]$ for any $c < 1.$ This can be viewed as partial MLR property when $\a > 1/2.$}
\end{Rq}

\end{document}